\documentclass[11pt,a4paper]{article}
\usepackage{amssymb,comment,amsmath,amsthm}
\usepackage{graphicx}
\usepackage{color}
\usepackage{titling,titlesec}
\usepackage{url}
\usepackage{setspace,amsthm,lipsum}

\newcounter{foo}

\newfont{\blb}{msbm10 scaled\magstep1}
\newfont{\comp}{cmr12 scaled\magstep1}
\newfont{\compb}{cmr10 scaled\magstep2}
\newfont{\sbb}{cmssbx10 scaled\magstep3}
\newfont{\sbbb}{cmssbx10 scaled\magstep5}
\newfont{\sbs}{cmssbx10 scaled\magstep1}
\newtheorem{thm}{Theorem}

\newtheorem{con}[foo]{Conjecture}

\newtheorem{prob}[foo]{Conjecture}

\title{\vspace{-1.5in} Recent progress in Ramsey Theory}

\date{ }

\author{
Jacques Verstraete\thanks{Department of Mathematics, University of California, San Diego, CA, 92093-0112 USA.
			Email: jverstraete@ucsd.edu.
			Research supported by NSF award DMS-1952786. \\
E-mail: jacques@ucsd.edu}}

\begin{document}

\maketitle

\vspace{-0.3in}

\begin{abstract}
\setstretch{1.05}{The classical Ramsey numbers $r(s,t)$ denote the minimum $n$ such that
every red-blue coloring of the edges of the complete graph $K_n$ contains either a red clique of order $s$ or a blue clique of
order $t$. These quantities are the centerpiece of graph Ramsey Theory, and have been studied for almost a century.
The Erd\H{o}s-Szekeres Theorem (1935) shows that for each $s \geq 2$, $r(s,t) = O(t^{s - 1})$ as $t \rightarrow \infty$.
We introduce a new approach using pseudorandom graphs
which shows $r(4,t) = \Omega(t^3/(\log t)^4)$ as $t \rightarrow \infty$, answering an old conjecture of Erd\H{o}s, and
we illustrate how to apply this approach to many other Ramsey and related combinatorial problems.}
\end{abstract}

\setcounter{tocdepth}{2}
\tableofcontents

\section{Introduction}

For integers $s,t \geq 2$, the centerpiece of graph Ramsey Theorey is the classical {\em Ramsey number} $r(s,t)$, which denotes the minimum $n$ such that every $n$-vertex graph contains either a clique of order $s$ or an independent set of size $t$. These quantities are the cornerstone of Ramsey Theory, and have been studied for many decades;
their existence is proved by Ramsey's Theorem~\cite{Ramsey} -- see also the book by Graham, Rothschild, Solymosi and Spencer~\cite{GRSS}. The quantities $r(t,t)$ are sometimes referred to as {\em diagonal Ramsey Numbers}, whereas $r(s,t)$ for fixed $s$ and $t \rightarrow \infty$ are referred to as {\em off-diagonal Ramsey numbers}. In a major recent breakthrough, Campos, Griffiths, Morris, Sahasrabudhe~\cite{CGMS} proved the following theorem:

\begin{thm}
There exists $\epsilon > 0$ such that for large enough $t$,
\begin{equation}
r(t,t) \leq (4 - \epsilon)^t
\end{equation}
\end{thm}

Explicit constructions giving lower bounds for Ramsey numbers have also been studied extensively, with many new results since the influential paper of Frankl and Wilson~\cite{FranklWilson}.

\medskip

The focus of this paper is to outline a new approach to Ramsey theory via pseudorandom graphs, based on the works~\cite{CMMV,MattheusV,MubayiV}. These and related techniques have been used effectively for not only estimating
Ramsey numbers~\cite{CMMV}, but also to {\em hypergraph Ramsey numbers}~\cite{KostochkaMubayiV}, {\em multicolor Ramsey numbers}~\cite{AlonRodl,HeWigderson}, {\em Ramsey minimal graphs}~\cite{MattheusV2}, {\em Erd\H{o}s-Rogers functions}~\cite{GJ,JanzerSudakov,MubayiV3} and {\em coloring hypergraphs}~\cite{MubayiV2}. A highlight of
this approach is that it determines $r(4,t)$ up to a factor of order $\log^2 \! t$, as proved by Mattheus and the author~\cite{MattheusV}:

\begin{thm}\label{r4t}
There exist constants $c_1,c_2 > 0$ such that for all $t \geq 3$,
\begin{equation}\label{r4tbound}
c_1\frac{t^3}{\log^4 \! t} \; \; \leq \; \; r(4,t) \; \; \leq \; \; c_2 \frac{t^3}{\log^2 \! t}.
\end{equation}
\end{thm}

This solves a conjecture of Erd\H{o}s~\cite{ChungGraham}, and also gives almost tight bounds on some multicolor Ramsey numbers.

\subsection{Ramsey numbers}

Since the original bounds of Erd\H{o}s and Szekeres~\cite{ErdosSzekeres} in 1935, estimates for $r(s,t)$ have been the subject of extensive research; for instance the conjecture $r(s,t) = t^{s-1}/\mbox{polylog}(t)$ as $t \rightarrow \infty$ dates back to at least 1947 -- see Erd\H{o}s~\cite{Erdos} and Chung and Graham~\cite{ChungGraham}. The only off-diagonal Ramsey numbers $r(s,t)$ for $s \geq 3$ whose order of magnitude
is known is $r(3,t)$, as it was shown by Kim~\cite{Kim} that $r(3,t) = \Omega(t^2/\log t)$ as $t \rightarrow \infty$, matching previous upper bounds
by Ajtai, Koml\'{o}s and Szemer\'{e}di~\cite{AjtaiKomlosSzemeredi} and Shearer~\cite{Shearer1,Shearer2}, and improving earlier bounds of Spencer~\cite{Spencer1,Spencer2}. The current state of the art is due Fiz Pontiveros, Griffiths and Morris~\cite{FGM} and Bohman and Keevash~\cite{BohmanKeevash1}, namely:
\begin{equation} \label{r3t}
(1 - o(1)) \frac{t^2}{4\log t} \; \;  \leq \; \; r(3,t) \; \; \leq \; \; (1 + o(1))\frac{t^{2}}{\log t}.
\end{equation}
Erd\H{o}s~\cite{ChungGraham} proposed the following conjecture:

\begin{con}\label{conjrst}
For all $s \geq 4$, there exists $c = c(s) > 0$ such that
\begin{equation}
r(s,t) \; \; \geq \; \; \frac{t^{s - 1}}{(\log t)^{c}}.
\end{equation}
\end{con}

The previous best bounds for $r(s,t)$ with fixed $s \geq 4$ are
\begin{equation} \label{rsn}
a \frac{t^{\frac{s+1}{2}}}{(\log t)^{\frac{s + 1}{2} - \frac{1}{s-2}}} \; \;  \leq \; \; r(s,t) \; \; \leq \; \; (1 + o(1)) \frac{t^{s-1}}{(\log t)^{s - 2}}
\end{equation}
where $a > 0$ is a constant depending only on $s$.
The lower bound comes from the so-called {\em $K_s$-free random process}, studied by Bohman and Keevash~\cite{BohmanKeevash2},
improving earlier bounds of Spencer~\cite{Spencer2}, while the upper bound comes from results of Li, Rousseau and Zang~\cite{LiRousseauZang}.
In light of Theorem \ref{r4t}, which solves Conjecture \ref{conjrst} for $s = 4$, we propose the following conjecture:

\begin{con}
The order of magnitude of $r(4,t)$ is $t^3/(\log t)^2$ as $t \rightarrow \infty$.
\end{con}

In the next section, we present general theorems connecting pseudorandom graphs, spectral graph theory and Ramsey theory,
which leads to Theorem \ref{r4t} as well as more general results on related Ramsey problems.

\subsection{Spectrally extremal graphs}

An {\em $(n,d,\lambda)$-graph} is an $n$-vertex $d$-regular graph whose adjacency matrix has eigenvalues $\lambda_1 \geq \lambda_2 \geq \dots \geq \lambda_n$ such that
$\max\{|\lambda_2|,|\lambda_n|\} = \lambda$. A survey of these graphs and related notions of pseudorandom graphs is given by Krivelevich and Sudakov~\cite{KrivelevichSudakov}. In the present context, the following theorem due to Mubayi and the author~\cite{MubayiV} connects $(n,d,\lambda)$-graphs with Ramsey Theory:

\begin{thm}\label{pseudo}
Let $F$ be a graph. If there exists an $F$-free $(n,d,\lambda)$-graph with $4e^2\lambda \geq (\log n)^2$, then
\begin{equation}\label{rftgeneral}
 r(F,t) \; \; > \; \; \frac{n(\log n)^2}{8e^2\lambda}.
 \end{equation}
\end{thm}

The key to effectively applying this theorem is to find an $F$-free $(n,d,\lambda)$-graph where $\lambda$ is as small as possible relative to $d$,
while $d$ is as large as possible relative to $n$. In general, the Alon-Boppana Theorem~\cite{Nilli1} shows that $\lambda = \Omega(\sqrt{d})$ when $d$ is not too close to $n$, so we refer to an $(n,d,\lambda)$-graph with $\lambda = O(\sqrt{d})$ as {\em spectrally extremal}.
Using the approach of Sudakov, Szab\'{o} and Vu~\cite{SudakovSzaboVu}, one can in theory then determine an upper bound on the value of $d$ relative to $n$ in an $F$-free
$(n,d,\lambda)$-graph; for example if $F = K_s$ then $d = O(n^{1 - 1/(2s - 3)})$. Remarkably, directly using Theorem \ref{pseudo}, this leads to the following theorem giving almost tight asymtotic bounds on $r(s,t)$ together with (\ref{rsn}):

\begin{thm}\label{rst}
If there exists a spectrally extremal $K_s$-free $(n,d,\lambda)$-graph with $d = \Omega(n^{1 - 1/(2s - 3)})$, then
\[ r(s,t) \; \; = \; \; \Omega\Bigl(\frac{t^{s - 1}}{(\log t)^{2s - 4}}\Bigr).\]
\end{thm}

This would prove Conjecture \ref{conjrst}. Such spectrally extremal triangle-free $(n,d,\lambda)$-graphs were constructed by Alon~\cite{AlonChung} and later by Kopparty~\cite{Kopparty}, but it is not known for any $s \geq 4$
whether there exists a spectrally extremal $K_s$-free graph with $d = n^{1 - 1/(2s - 3) - o(1)}$. The current record is
due to Bishnoi, Ihringer and Pepe~\cite{BishnoiIhringerPepe}, with $d = \Omega(n^{1 - 1/(s - 1)})$ for $s \geq 4$.
Since there exist spectrally extremal triangle-free $n$-vertex graphs with $d = \Theta(n^{2/3})$, as constructed in~\cite{Alon,Kopparty},
an application of Theorem \ref{rst} gives $r(3,t) = \Omega(t^2/(\log t)^2)$, which matches the lower bound obtained by direct application of
the Lov\'{a}sz local lemma~\cite{Spencer1}. The best explicit constructions giving lower bounds on $r(s,t)$ for $s \in \{4,5,6\}$ are described by
Kostochka, Pudlak and R\"{o}dl~\cite{KostochkaPudlakRodl}.

\medskip

More generally, for a graph $F$ and an integer $t \geq 2$, let $r(F,t)$ denote the minimum $n$ such that every $n$-vertex $F$-free graph contains an independent set of size $t$. The notorious {\em cycle complete-graph Ramsey numbers} are concerned with the case $F = C_k$, that is, $F$ is a cycle of length $k$. Perhaps one of the most challenging open conjectures in Ramsey Theory is the following conjecture of Erd\H{o}s~\cite{ChungGraham}:

\begin{con}
There exists $\varepsilon > 0$ such that if $t$ is large enough, then
\begin{equation}\label{rc4t}
r(C_4,t) \; \; \leq \; \; t^{2 - \varepsilon}.
\end{equation}
\end{con}

The current best bounds for $r(C_4,t)$ are $r(C_4,t) = O(t^2/(\log t)^2)$ due to Szemer\'{e}di (see~\cite{EFRS}) and $r(C_4,t) = \Omega(t^{3/2}/\log t)$, the latter coming from the random $C_4$-free process. A {\em polarity graph} of a projective plane of order $q$ may be viewed as an $(n,d,\lambda)$-graph where
$n = q^2 + q + 1$, $d = q + 1$ and $\lambda = \sqrt{q}$, and these graphs do not contain cycles of length four. Applying Theorem \ref{pseudo} to these graphs,
one immediately obtains $r(C_4,t) = \Omega(t^{3/2}/\log t)$. Using finite geometry, we obtain in a similar way the following bounds on $r(C_{\ell},t)$, which exceed the previous best known bounds from the random $C_{\ell}$-free process:

\begin{thm} \label{cycles}
As $t \rightarrow \infty$, 	
\begin{equation}\label{evencycles}
r(C_6, t) \; \; = \; \;  \Omega\left(\frac{t^{5/4}}{\log^{1/2} t}\right) \qquad \hbox{ and } \qquad r(C_{10}, t) \; \; = \; \;  \Omega\left(\frac{t^{9/8}}{\log^{1/4} t}\right).
\end{equation}	
\end{thm}

These results are obtained by considering polarity graphs~\cite{LazebnikUstimenkoWoldar}  of {\em generalized quadrangles} and {\em generalized hexagons}.
 The best upper bounds on cycle-complete graph Ramsey numbers
are due to Caro, Li, Rousseau and Zhang~\cite{CLRZ} for even cycles and Sudakov~\cite{Sudakov} for odd cycles.
For a survey on extremal problems for cycles in graphs, see~\cite{V}.

\subsection{Transference from extremal to Ramsey}

 The lack of dense spectrally extremal $F$-free graphs, and in particular $K_4$-free graphs, precludes an application of Theorem \ref{pseudo} to obtain Theorem \ref{r4t}.
Instead, when $F$ is a non-bipartite graph, we use a weaker notion of pseudorandom graphs combined with the following approach in the language of hypergraphs~\cite{CMMV}.
Let $r \geq 3$ and let $H$ be an $r$-uniform hypergraph in which no two edges intersect in at least two vertices. Let $\partial H$ denote the graph whose vertex set is $V(H)$ and whose edge set consists of pairs of vertices in a common edge of $H$ -- this is the {\em shadow graph} of $H$. We say that $H$ is {\em strongly $F$-free} if for every isomorphic copy $F' \subseteq \partial H$ of $F$, there exists an edge $e \in E(H)$ such that the subgraph of $F'$ induced by $e \cap V(F')$ is non-bipartite. For instance, if $F$ is a triangle then $H$ does not contain three edges $e_1,e_2,e_3$ such that $|e_i \cap e_j| = 1$ for $i \neq j$ and $e_1 \cap e_2 \cap e_3 = \emptyset$. The following theorem~\cite{CMMV} connects the problem of finding extremal strongly $F$-free hypergraphs to bounds on Ramsey numbers:

\begin{thm}\label{transfer}
Let $F$ be a non-bipartite graph, and suppose there exists an $r$-uniform $n$-vertex $d$-regular strongly $F$-free hypergraph, where $r \geq 2048(\log n)^3$.
If $t = \lceil 256n(\log n)^2/rd \rceil$, then
\begin{equation}
r(F,t) \; \; = \; \; \Omega\Bigl(\frac{rt}{\log n}\Bigr).
\end{equation}
\end{thm}

We shall show how Theorem \ref{r4t} may be derived from Theorem \ref{transfer}.  The approach used to prove Theorem \ref{r4t} gives the current best bounds~\cite{CMMV,MubayiV} on cycle-complete graph Ramsey numbers $r(C_k,t)$ for all odd values of $k$ and $k \in \{4,6,10\}$. For instance, in the case of the pentagon,
we have the following result~\cite{CMMV} -- the upper bound in the following theorem is due to Sudakov~\cite{Sudakov}:

\begin{thm}\label{rc5t}
There exist constants $c_1,c_2 > 0$ such that for all $t \geq 3$,
\[ c_1 \frac{t^{10/7}}{(\log t)^{13/7}} \; \; \leq \; \; r(C_5,t) \; \; \leq \; \; \frac{t^{3/2}}{(\log t)^{1/2}}.\]
\end{thm}

The proof of Theorem \ref{rc5t} uses finite geometric object known as a twisted triality hexagon~\cite{CMMV},
which in the language of hypergraphs allows an application of Theorem \ref{transfer}. This theorem suggests that it is plausible that $r(C_5,t)$ has order of magnitude close to $t^{3/2}$:

\begin{con}
$r(C_5,t)  = t^{3/2 - o(1)}$ as $t \rightarrow \infty$.
\end{con}

\section{Pseudorandom graphs}

If $A$ is a square symmetric matrix, then the eigenvalues of $A$ are real. When $A$ is an $n$ by $n$ matrix, we denote them $\lambda_1 \geq \lambda_2 \geq \dots \geq \lambda_n$. If the corresponding eigenvectors forming an orthonormal basis are $e_1,e_2,\dots,e_n$, respectively, then for any $x \in \mathbb R^n$, we may write $x = \sum x_i e_i$ and
\begin{equation}\label{easy} \left<Ax,x\right> = \sum_{i = 1}^n \lambda_i x_i^2.
\end{equation}
When $A$ is the adjacency matrix of a graph $G$, we define $\lambda_i(G) = \lambda_i$ and
\begin{equation}\label{lambdadef}
\lambda(G) = \max\{|\lambda_i| : 2 \leq i \leq n\}
\end{equation}
and we refer to the eigenvalues of the graph rather than the eigenvalues of $A$.
For $X \subseteq V(G)$, let $e(X)$ denote the number of edges $\{u,v\} \in E(G)$ with $u,v \in X$. If $x$ is the characteristic vector of $X$
and $x = \sum x_i e_i$, then
\begin{equation}\label{edgeadj}
2e(X) = \left<Ax,x\right>.
\end{equation}
If $A$ is doubly stochastic, then $e_1$ is the constant unit vector with eigenvalue equal to the common row sum.
In particular, if $A$ is the adjacency matrix of a $d$-regular graph $G$, then $\lambda_1(G) = d$ and if $\lambda(G) = \lambda$,
then $G$ is an $(n,d,\lambda)$-graph. We shall see shortly that if all the eigenvalues of $G$ except $\lambda_1(G)$ are small in absolute value,
then large subsets of the vertices induce close to the expected number of edges -- this is the expander mixing lemma -- and hence such graphs could be referred to collectively as pseudorandom graphs~\cite{KrivelevichSudakov}.

\subsection{The Alon-Boppana Theorem}

The infinite $d$-regular tree $T_d$ is the universal
cover of $d$-regular graphs, so for any $d$-regular graph $G$ with $n$ vertices, the number of closed walks of length $2k$ from $x \in V(G)$ to $x$ is at least the number of closed walks of length $2k$ from the root of $T_d$ to the root of $T_d$. This is at least
\begin{equation}\label{ab-trees}
 \frac{1}{k}{2k - 2 \choose k - 1} d(d - 1)^{k - 1}.
 \end{equation}
The total number of walks of length $2k$ in $G$ is at most $d^{2k} + (n - 1)\lambda^{2k}$ by the trace formula, so
\begin{equation}\label{ab-inequality}
 d^{2k} + (n - 1)\lambda^{2k} \geq \frac{n}{k}{2k - 2 \choose k - 1} d(d - 1)^{k - 1}.
 \end{equation}
For fixed $d$, by selecting $k$ to depend appropriately on $n$ in (\ref{ab-inequality}), we obtain the Alon-Boppana Theorem~\cite{Nilli1}:

\begin{thm} Let $d \geq 1$. If $G_n$ is a $d$-regular $n$-vertex graph then
\begin{equation}\liminf_{n \rightarrow \infty} \lambda(G_n) \geq 2\sqrt{d - 1}.\end{equation}
\end{thm}

A very short proof of the Alon-Boppana Theorem was found by Alon~\cite{Nilli1,Nilli2}. In light of this theorem,
we refer to an $(n,d,\lambda)$ graph with $\lambda = O(\sqrt{d})$ as a {\em spectrally extremal graph}.
For instance, {\em Ramanujan graphs} are spectrally extremal $d$-regular graphs with $\lambda \leq 2\sqrt{d - 1}$, and were constructed by Lubotzky, Philips and Sarnak~\cite{LubotzkyPhillipsSarnak} and Margulis~\cite{Margulis} as Cayley graphs for certain prime powers $d$ (bipartite constructions with $\max\{|\lambda_2|,|\lambda_{n-1}|\} \leq 2\sqrt{d - 1}$ were given for all values of $d$ via polynomial interlacing by Marcus, Spielman and Srivastava~\cite{MarcusSpielmanSrivastava}). It turns out that random $d$-regular graphs are also ``asymptotically'' Ramanujan graphs with high probability, as proved in a major work by Friedman~\cite{Friedman}.

\subsection{The Expander Mixing Lemma}

Let $G$ be an $(n,d,\lambda)$-graph. If $X \subseteq V(G)$ has characteristic vector $x$ and $x = \sum_{i = 1}^n x_i e_i$
with $\{e_1,e_2,\dots,e_n\}$ an orthonormal basis of eigenvectors, then
\begin{equation}\label{expand}
\left<Ax,x\right> = \sum_{i = 1}^n \lambda_i x_i^2.
\end{equation}
Since $e_1$ is the constant vector, $x_1 = \left<x,e_1\right> = |X|/\sqrt{n}$. Estimating the remaining terms in (\ref{expand})
using $|\lambda_i| \leq \lambda$, and applying (\ref{edgeadj}), we obtain the {\em expander mixing lemma} of Alon~\cite{AlonChung}:

\begin{thm}
Let $G$ be an $(n,d,\lambda)$-graph and let $X \subseteq V(G)$. Then
\begin{eqnarray}
\lambda_n |X| \leq 2e(X)  - \tfrac{d}{n}|X|^2 \leq \lambda |X|
\end{eqnarray}
\end{thm}

In particular, the smaller the value of $\lambda$, the more tightly concentrated the number of edges in subsets $X$ is around $d|X|^2/2n$.

\subsection{Spectrally extremal graphs}

Mantel's Theorem~\cite{Mantel} states that a triangle-free graph with $2n$ vertices has at most $n^2$ edges, with equality only for complete $n$ by $n$ bipartite
graphs, which are triangle-free $(2n,n,n)$-graphs.  In general, one may ask for triangle-free $(n,d,\lambda)$-graphs, where $d,n$ and $\lambda$ must satisfy
(\ref{ab-inequality}), and in particular, $\lambda = \Omega(\sqrt{d})$. This leads to an extremal problem: when $\lambda = O(\sqrt{d})$,
what is the largest $d$ for which there is a triangle-free $(n,d,\lambda)$-graph, or more generally an $F$-free $(n,d,\lambda)$-graph? If $G$ is any triangle-free $(n,d,\lambda)$-graph with adjacency matrix $A$,  then
\begin{equation}\label{hoffman}
0 = \hbox{tr}(A^3) \ge  d^3 - \lambda^3(n-1)
\end{equation}
and so $d \leq \lambda (n - 1)^{1/3}$. If $\lambda = O(\sqrt{d})$, then this gives $d = O(n^{2/3})$. Alon~\cite{Alon}
gave an ingenious construction of a triangle-free Cayley $(n,d,\lambda)$-graph with $d = \Omega(n^{2/3})$ and $\lambda = O(n^{1/3})$, showing
the the above bounds are tight (see also Conlon~\cite{Conlon} and Kopparty~\cite{Kopparty}):

\begin{thm}\label{alon-triangle}
There exist triangle-free $(n,d,\lambda)$-graphs with $\lambda = O(n^{1/3})$ and $d = \Omega(n^{2/3})$ as $n \rightarrow \infty$.
\end{thm}

Similar arguments show that if $k$ is odd and $G$ is a $C_{k}$-free graph then $d = O(\lambda n^{1/k})$ and when $\lambda = O(\sqrt{d})$ this gives $d = O(n^{2/k})$, and such spectrally extremal graphs exist~\cite{AlonKahale}. If $F$ is a bipartite graph, then
in many cases, the densest known $n$-vertex $F$-free graphs are graphs whose degrees are all close to some number $d$ and with $\lambda = O(\sqrt{d})$ -- see F\"{u}redi and Simonovits~\cite{FurediSimonovits} for a survey of bipartite extremal problems. A particular rich source of such graphs is the random polynomial model, due to Bukh and Conlon~\cite{BukhConlon}, as well as {\em strongly regular graphs} -- see Brouwer and van Maldeghem~\cite{BrouwerVanMaldeghem} for a survey.
Another fruitful source of examples is {\em Cayley graphs}, since the eigenvalues are determined by character sums.

\subsection{Spectrally extremal clique-free graphs}

More generally, Sudakov, Szabo and Vu~\cite{SudakovSzaboVu} showed
\begin{equation} d = O(\lambda^{\frac{1}{s  - 1}} n^{1 - \frac{1}{s  - 1}}).
\end{equation}
when $G$ is a $K_s$-free $(n,d,\lambda)$-graph. When $\lambda = O(\sqrt{d})$ this gives
\begin{equation}
d = O(n^{1 - \frac{1}{2s - 3}}).
\end{equation}
 Constructing such spectrally extremal $K_s$-free graphs for $s \geq 4$ is considered to be one of the major open problems in the area.
 Examples of constructions were given by Alon and Krivelevich~\cite{AlonKrivelevich} and Krivelevich~\cite{Krivelevich}, and the current best constructions were given by Bishnoi, Ihringer and Pepe~\cite{BishnoiIhringerPepe},
who define a graph $\Gamma[q,s]$ using the bilinear form
$Q(x,y) = \xi x_1 y_1 + \sum_{i = 2}^{s+1} x_i y_i$ with $\xi$ a quadratic non-residue in $\mathbb F_q$. Let $\chi$ be the quadratic character of $\mathbb F_q$.
The vertex set of $\Gamma[q,s]$ is a set of one-dimensional subspaces of $\mathbb F_q^{s + 1}$, defined as follows:
\begin{equation}
\{\left<x\right> < \mathbb F_q^{s+1} : \chi(Q(x,x)) = 1\}
\end{equation}
and two one-dimensional subspaces $\left<x\right>$ and $\left<y\right>$ are adjacent if $\chi(Q(x,y)) = 1$. A key fact is that the neighborhood of a
vertex in $\Gamma[q,s]$ induces $\Gamma[q,s-1]$, and $\Gamma[q,1]$ is empty with $q/2$ vertices. Then an induction shows $\Gamma[q,s]$ is $K_{s + 1}$-free,
and $\Gamma[q,s]$ is an $(n,d,\lambda)$-graph with $d = \Theta(n^{1 - \frac{1}{s}})$ and $\lambda = O(\sqrt{d})$.

\subsection{Independent sets}

An {\em independent set} in a graph $G$ is simply a set $X$ of vertices such that no edge of $G$ is contained in $X$.
From the expander mixing lemma, one may obtain the {\em Hoffman-Delsarte bound} or {\em ratio bound}, which states that
if $X$ is an independent set in a $d$-regular graph $G$ whose adjacency matrix has eigenvalues $d = \lambda_1 \geq \lambda_2 \geq \dots \geq \lambda_n$ and 
$\lambda = \max\{|\lambda_i| : 2 \leq i \leq n\}$, then  $|X| \leq -n\lambda_n/(d - \lambda_n) \leq \lambda n/(d + \lambda)$ -- see Haemers~\cite{Haemers}. 
We shall see shortly that the number of independent sets of size not much larger than $n/(d + 1)$ in such graphs is very small.
It turns out that this counting of independent sets in {\em pseudoramdom graphs} is key to Ramsey theoretic arguments. For our purposes,
an $n$-vertex graph $G$ is {\em $(\alpha,m)$-pseudorandom} if for every set $X \subseteq V(G)$ of size at least $m$,
\begin{equation}\label{condition}
e(X) \geq \alpha {|X| \choose 2}.
\end{equation}
According to the expander mixing lemma, an $(n,d,\lambda)$-graph is indeed $(\alpha,m)$-pseudorandom
for some $m = \Theta(n\lambda/d)$ and $\alpha = \Theta(d/n)$ as $d \rightarrow \infty$. By Tur\'{a}n's Theorem~\cite{Turan}, every graph of average degree $d$ has an independent set of size at least $t = n/(d + 1)$. Furthermore, the number of such independent sets in any $n$-vertex $d$-regular graph (when $d + 1 | n$) is at least
\begin{equation}
\frac{1}{t!} n(n - d - 1)(n - 2d - 2) \dots = (d + 1)^t.
\end{equation}
A disjoint union of cliques $K_{d + 1}$ shows this is tight.
In pseudorandom graphs, it turns out there are few independent sets that are much larger than this.
The following theorem due to Kohayakawa, Lee, R\"{o}dl and Samotij~\cite{KLRS}, which is
an element of the {\em method of containers}~\cite{BaloghMorrisSamotij,SaxtonThomason},
implies a similar earlier theorem of Alon and R\"{o}dl~\cite{AlonRodl} for $(n,d,\lambda)$-graphs:

\begin{thm}\label{alonrodl}
Let $n \geq m \geq t \geq s \geq 1$ and $\alpha \in [0,1]$, and let $G$ be an $n$-vertex $(\alpha,m)$-pseudorandom graph, and 
$\exp(-\alpha s)n \leq m$. Then the number of independent sets of size $t \geq s$ in $G$ is at most
\begin{equation} \label{finalcount}
{n \choose s} {m \choose t - s}.
\end{equation}
In particular, if $G$ is an $(n,d,\lambda)$-graph, then the number of independent
sets of size $t \geq 2n(\log n)^2/d$ is at most
\begin{equation} \label{finalcount2}
\Bigl(\frac{4e^2 \lambda}{\log^2\! n}\Bigr)^t.
\end{equation}
\end{thm}

The statement (\ref{finalcount2}) follows from (\ref{finalcount}) using the fact that
an $(n,d,\lambda)$-graph is $(\alpha,m)$-pseudorandom for some $m = \Theta(n\lambda/d)$ and $\alpha = \Theta(d/n)$ as $d \rightarrow \infty$,
and choosing roughly $r = \Theta(t/\log n)$. Alon and R\"{o}dl~\cite{AlonRodl} used their theorem to obtain remarkably good bounds on
multicolor Ramsey numbers, by simply randomly superimposing certain $(n,d,\lambda)$-graphs on top of each other.
It is not clear whether the above theorem is sharp up to constant factors, for instance we may ask the following:

\begin{prob}
Let $G$ be an $(n,d,\lambda)$-graph. Does there exist a constant $C > 0$ such that the number of independent
sets of size $t \geq C(n\log n)/d$ in $G$ is at most $(C\lambda/\log n)^t$?
\end{prob}

A positive answer to this question would for instance give $r(3,t) = \Omega(t^2/\log t)$,
and under the spectral conditions of Theorem \ref{rst}, one would obtain $r(s,t) = \Theta(t^{s - 1}/(\log t)^{s - 2})$.
A survey on counting independent sets is given by Samotij~\cite{Samotij}.

\section{Proofs of Theorems}

\subsection{Proof of Theorem \ref{alonrodl}}

We aim to prove (\ref{finalcount}). Let $N_G[v] = N_G(v) \cup \{v\}$ for $v \in V(G)$. 
The main claim is that for each independent set $I \subseteq V(G)$ of size $t$, there exists $i \leq s$ and vertices $v_1,v_2,\dots,v_i \in I$ such that
\[ \Bigl|V(G) \;  \setminus  \; \bigcup_{j \leq i} N_G[v_j]\Bigr| \leq m.\]
Let $v_1$ be the first vertex of $I$ in a max degree ordering $>_0$ of $G_0 = G$. Let $n_0 = n$, and 
\[ n_1 = |\{v \in V(G_0) : v_1 >_0 v\}| \qquad \mbox{ and } \qquad G_1 = G_0 - N_{G_0}[v_1].\] 
Suppose we have defined vertices $v_1,v_2,\dots,v_j \in I$, and $G_j = G_{j - 1} - N_{G_{j - 1}}[v_j]$. Then let $v_{j + 1}$ be the first 
vertex of $I$ in the max degree ordering $>_j$ of $G_j$, and let 
\[ n_{j + 1} = |\{v \in V(G_i) : v_{j+1} >_j v\}| \qquad \mbox{ and } \qquad G_{j + 1} = G_j - N_{G_j}[v_{j+1}].\]
The claim is proved if $n_i \leq m$ for some $i \in [0,s]$, so we assume $n_i > m$ for all $i \in [0,s]$. 
Let $X_i = \{v \in V(G_{i - 1}) \backslash \{v_i\} : v_i >_{i - 1}v\}$. Then $|X_i| = n_i - 1 \geq m$, so it follows that for all $i \in [0,s]$, 
\[ |N_{G_{i - 1}}[v_i]| \geq \max\{|N_{G_{i - 1}}(v)| : v \in X_i\} + 1 \geq \alpha (|X_i| - 1) + 1 \geq  \alpha n_i.\] 
In particular, $n_{i + 1} \leq (1 - \alpha)n_i$ for all $i \in [0,s-1]$. Using $\exp(-\alpha s)n \leq m$, this implies  
\[ m < n_s \leq (1 - \alpha)n_{s-1} \leq \dots \leq (1 - \alpha)^{s - 1} n_1 \leq (1 - \alpha)^s n_0 \leq \exp(-\alpha s)n \leq m.\]
This contradiction proves the claim. \qed

\subsection{Proof of Theorem \ref{pseudo}}

We are to show that if there exists an $F$-free $(n,d,\lambda)$-graph $G$, then when $t \geq  2n(\log n)^2/d$,
\[ r(F,t) > \frac{n(\log n)^2}{8e^2\lambda}.\]
By Theorem \ref{alonrodl}, the number of independent sets of
size $t \geq 2n(\log d)^2/d$ in $G$ is at most
\[ \Bigl(\frac{4e^2 \lambda}{\log^2\! n}\Bigr)^t.\]
Sample vertices of $G$ independently with probability $p = (\log^2\! n)/(4e^2\lambda)$ and let $X$ be the set of sampled vertices.
Note that since $\lambda \geq (\log^2\! n)/4e^2$, this is valid.
Then $\mbox{E}(X) = pn$ and the expected number of independent sets of size $t$ in $X$ is $1$. Removing one vertex from each independent
set of size $t$ in $X$ gives a set $X'$ of $pn - 1$ vertices with no independent set of size $t$ and inducing no copy of $F$. Therefore
\[ r(F,t) > pn - 1 > \frac{n(\log n)^2}{8e^2\lambda}.\]
This proves Theorem \ref{pseudo}.

\subsection{Proof of Theorem \ref{rst}}

To derive Theorem \ref{rst}, simply apply Theorem \ref{pseudo} to a spectrally extremal $n$-vertex $d$-regular $K_s$-free graph,
where $d = \Theta(n^{1 - 1/(2s - 3)})$.

\subsection{Proof of Theorem \ref{transfer}}

The proof of Theorem \ref{transfer} involves constructing a graph based on the shadow graph of a strongly $F$-free hypergraph $H$ as stated in the theorem. We observe
that for any set $X \subseteq V(G)$,
\[ \sum_{e \in E(H)} |e \cap X| = \sum_{x \in X} d(x)\]
where $d(x)$ is the number of edges in $H$ containing $x$. Since $d(x) \geq d$ for all $x \in V(H)$, we find that on average,
for each edge of $H$, $|e \cap X| \geq r|X|/n$. If $|X| \geq 2n/r$, this implies that the number of edges of $\partial H$ contained in $X$
is
\[ \sum_{e \in E(H)} {|e \cap X| \choose 2} \geq \frac{dn}{r} {\frac{r|X|}{n} \choose 2} \geq \frac{dr}{2n} {|X| \choose 2}.\]
In particular, the shadow graph is $(\alpha,m)$-pseudorandom with $\alpha = dr/2n$ and $m = 2n/r$.
The next step is to define a random $F$-free subgraph of the shadow graph. Independently for each edge of $H$, randomly color the vertices of that edge with two colors, and let $G$ be obtained from the shadow graph $\partial H$ by taking all edges which have ends of different colors. Since $H$ is strongly $F$-free,
we conclude that $G$ is $F$-free. The main task is to show that $G$ inherits the pseudorandomness properties of $\partial H$. Indeed, we expect
about half the edges of $\partial H$ to be edges of $G$ in any subset $X$. The details, which are given in~\cite{CMMV}, involve the use
of Hoeffding's martingale inequality~\cite{Hoeffding}, and the conclusion is that with high probability, $G$ is $(\alpha,m)$-pseudorandom with $\alpha = \Theta(dr/n)$ and
$m = \Theta(n/r)$ as $n \rightarrow \infty$. According to Theorem \ref{alonrodl}, specifically (\ref{finalcount}),
there are at most ${n \choose s}{m \choose t - s}$ independent
sets of size $t \geq s$ in $G$, when $\exp(-\alpha s)n \leq m$.  We choose $t$ as given in the theorem and $s = \Theta(t/\log n)$ so that
the number of independent sets of size $t$ is at most $\Theta(n\log n/rt)^t$. Sampling vertices independently with suitable probability
$\Theta(rt/n\log n)$, we destroy all these independent sets to obtain an $F$-free induced subgraph with $\Omega(rt/\log n)$ vertices and no independent
set of size $t$, thereby showing $r(F,t) = \Omega(rt/\log n)$.

\subsection{Proof of Theorem \ref{r4t}}

Theorem \ref{r4t} is proved using Theorem \ref{transfer}. The hypergraph to which we apply Theorem \ref{transfer} to obtain the claimed bounds on $r(4,t)$ has vertex set equal to the set of lines of the {\em Hermitian unital} and
edges consisting of all lines through a point of the unital: this is a $K_4$-free $n$-vertex $r$-uniform $d$-regular hypergraph with $n$ vertices, where for some prime power $q$, we have $r = q^2$,  $n = q^2(q^2 - q + 1)$, $d = q + 1$. For background on unitals, see Barwick and Ebert~\cite{BarwickEbert}. In particular, with $t = \lceil 256n (\log n)^2/rd \rceil = \Theta(q(\log q)^2)$, we have
\[ r(4,t) = \Omega\Bigl(\frac{rt}{\log^2 n}\Bigr) = \Omega\Bigl(\frac{t^3}{(\log t)^4}\Bigr).\]

\subsection{Proof of Theorem \ref{rc5t}}

We apply Theorem \ref{transfer} similarly to obtain Theorem \ref{rc5t}: here $H$ is the strongly $C_5$-free hypergraph whose vertices are the lines of
the {\em twisted triality hexagon} and whose edges are the set of lines through each point of the hexagon. The parameters for prime powers $q$ are $n = (q^3 + 1)(q^8 + q^4 + 1)$, $r = q^3 + 1$ and $d = q + 1$. Then $H$ is $C_5$-free and so Theorem \ref{transfer} with $t = \Theta(q(\log q)^2)$ gives
\[ r(C_5,t) = \Omega\Bigl(\frac{(q^3 + 1)t}{\log^2 n}\Bigr) = \Omega\Bigl(\frac{t^{10/7}}{(\log t)^{13/7}}\Bigr).\]
This supersedes the $C_5$-free process by a polynomial in $t$.

\section{Concluding remarks}

The techniques presented in this paper have been used effectively for not only estimating
Ramsey numbers~\cite{CMMV}, but also to {\em hypergraph Ramsey numbers}~\cite{KostochkaMubayiV}, {\it multicolor Ramsey numbers}~\cite{AlonRodl,HeWigderson}, {\em Erd\H{o}s-Rogers functions}~\cite{GJ,JanzerSudakov,MubayiV}, {\em degree minimal Ramsey graphs}~\cite{MattheusV2} and {\em coloring hypergraphs}~\cite{MubayiV3}. In this section, we give a brief description of new results.

\medskip

$\bullet$ Let $T$ denote the hypergraph consisting of three triples $e,f,g$ such that $e \cap f \cap g = \emptyset$ whereas $e \cap f,f \cap g$ and $e \cap g$ are all non-empty.
This is sometimes referred to as a {\em generalized triangle}. Let $r(T,t)$ denote the minimum $n$ such that every $T$-free triple system on $n$ vertices contains
a set of $t$ vertices containing no triples -- an {\em independent set}. By randomly placing a suitable family of triples inside each line of a generalized quadrangle,
the following result is proved in~\cite{KostochkaMubayiV}:

\begin{thm}
There exist absolute constants $c_1,c_2 > 0$ such that for all $t \geq 3$,
\[ c_1 \frac{t^{3/2}}{(\log t)^{3/4}} \; \; \leq \; \; r(T,t) \; \; \leq \; \;  c_2 t^{3/2}.\]
\end{thm}

The main open problem is to prove $r(T,t) = o(t^{3/2})$, as conjectured in~\cite{KostochkaMubayiV}.

\medskip

$\bullet$ Using Theorem \ref{r4t}, one can obtain almost tight lower bounds on $k$-color Ramsey numbers: let $r_k(4,t)$ denote the minimum $n$ such that
in every $k$-coloring of the edges of $K_n$, one has either a monochromatic $K_4$ in one of the first $k - 1$ colors or a monochromatic $K_t$ in the last color.
Note that $r_2(4,t) = r(4,t)$. Then Theorem \ref{r4t} and the upper bounds of He and Wigderson~\cite{HeWigderson} give the following theorem:

\begin{thm}
There exist absolute constants $c_1,c_2 > 0$ such that for all $t \geq 3$,
\[ c_1 \frac{t^{2k - 1}}{(\log t)^{4k-4}} \; \; \leq \; \; r_k(4,t) \; \; \leq \; \;  c_2 \frac{t^{2k - 1}}{(\log t)^{2k-2}}.\]
\end{thm}

\medskip

$\bullet$ The construction in~\cite{MattheusV} has been used to derive bounds on quantities known as {\em Erd\H{o}s-Rogers functions}. For an integer $t > s \geq 2$, the {\em Erd\H{o}s-Rogers function} $f_{s,t}(n)$ is the
maximum integer $m$ such that every $n$-vertex $K_{t}$-free graph has a $K_s$-free induced subgraph with
$m$ vertices. These quantities were introduced by Erd\H{o}s and Rogers~\cite{ER} more than sixty years ago. These quantities are generalizations of Ramsey numbers,
since $r(s,t)$ and $f_{2,t}(n)$ determine one another for each $t \geq 2$.
New upper bounds on $f_{s,s+2}(n)$ were recently obtained by Janzer and Sudakov~\cite{JanzerSudakov} using the construction in~\cite{MattheusV} and the method
of containers~\cite{BaloghMorrisSamotij,SaxtonThomason}, as well as new results by Gowers and Janzer~\cite{GJ} on $f_{s,t}(n)$. However, there is no pair $(s,t)$ with $t \geq s + 2 \geq 4$ for which it is known that
$f_{s,t}(n) = n^{\alpha + o(1)}$ for some $\alpha = \alpha(s,t)$. We focus here on $f_{s,s+1}(n)$. Since the quantities $r(3,t)$ are known to within a constant factor, one deduces $f_{2,3}(n)$ has order of magnitude $\sqrt{n\log n}$ as $n \rightarrow \infty$. Building on earlier works of Dudek and R\"{o}dl~\cite{DR}, Dudek, Retter and R\"{o}dl~\cite{DRR} and  Wolfovits~\cite{Wolfovits}, we significantly improve the upper bounds on $f_{s,s+1}(n)$ in~\cite{MubayiV3} as follows:

\begin{thm}\label{er}
For each fixed $s \geq 3$,
\[
f_{s,s+1}(n) \; \; \leq \; \; 2^{100s} \sqrt{n} \, \log n.
\]
\end{thm}

A lower bound $f_{s,s+1}(n) = \Omega(\sqrt{n\log n}/\log \log n)$ can be obtained from results of Shearer~\cite{Shearer2}.
Therefore this theorem incidentally shows $f_{s,s+1}(n) = n^{1/2 + o(1)}$ for $s = o(\log n)$. Dudek, Retter and R\"{o}dl~\cite{DRR} asked whether for $s \geq 3$,  \[ \lim_{n \rightarrow \infty} \frac{f_{s,s+1}(n)}{f_{s - 1,s}(n)} = \infty\]
as $n \rightarrow \infty$, which would not be the case for $s \geq 4$ if our upper bound is tight up to a constant factor. We tentatively make the following conjecture:

\begin{con}\label{one}
For all $s \geq 3$, as $n \rightarrow \infty$,
\[ f_{s,s+1}(n) = \sqrt{n} (\log n)^{1 - o(1)}.\]
\end{con}

This would be in contrast to $f_{2,3}(n)$, which has order of magnitude $\sqrt{n\log n}$ as $n \rightarrow \infty$. A key barrier in proving lower bounds
on say $f_{3,4}(n)$ is showing that in a suitably dense $K_4$-free graph one can typically find induced triangle-free subgraphs with substantially more vertices than the independence number
of the graph.

\medskip

$\bullet$ A graph $G$ is {\em $r$-Ramsey minimal} for a graph $H$ if every $r$-coloring of $E(G)$ contains a monochromatic copy of $H$, but no proper subgraph of $G$ possesses this property. The smallest minimum degree $s_r(H)$ over all $r$-Ramsey minimal graphs $G$ was first studied by Burr, Erd\H{o}s, and Lov\'{a}sz~\cite{BurrErdosLovasz}, who showed $s_2(K_k) = (k-1)^2$. Using the construction in~\cite{MattheusV}, we improve the most recent upper bounds on $s_r(K_k)$ due to Fox, Grinshpun, Liebenau, Person and
Szab\'{o}~\cite{FGLPS} and Bishnoi and LesGourges~\cite{BishnoiLesGourges} as follows:

\begin{thm}\label{ramseymin}
There exists an absolute constant $c > 0$ such that
for $k,r \geq 3$, $s_r(k) \leq 2^{100k} r^2 \log r$ and
\[ s_r(k) \leq (rk)^{2 + \frac{c}{k}}[(\log k)^2 + (\log k)^{1 + \frac{c}{k}}(\log r)^{1 + \frac{c}{k}}].\]
On the other hand, for all $k,r \geq 3$,
\[ s_r(k) = \Omega((k - 1)r^2).\]
\end{thm}

When $k \geq \log r$, this theorem gives $s_r(k) = O(r^2k^2(\log k)^2\log r)$, which improves earlier bounds
$s_r(k) = O(r^3 k^3(\log k)^3)$ of Bishnoi and LesGourges~\cite{BishnoiLesGourges}.
The main open problem is to find better lower bounds on $s_r(k)$, especially when $r$ and $k$ are both large -- the order of magnitude is
at most $(rk\log k)^2 \log r$ from the above theorem when $k \geq \log r$, whereas the lower bound is of order $(k - 1)r^2$.

\medskip

$\bullet$ Let $q \geq k \geq 2$. Consider a hypergraph $H$ which comes from taking all $k$-element subsets of sets in a family $F$ of $e$ sets of size $q$ which pairwise
intersect in at most one vertices.
Such hypergraphs arise naturally in finite geometry, for instance, if $F$ is the set of lines of a projective plane then $e = q^2 + q + 1$,
and for $k = 3$, then the hypergraph $H$ is the family of all collinear triples. In the case $k = 2$, the Erd\H{o}s-Faber-Lov\'{a}sz conjecture
states that if $|F| = q$ then the vertex chromatic number of the graph $H$ is $q$. In other words, the edge-disjoint union of $q$ cliques of size $q$ has chromatic number $q$. This was proved by Kang, Kelly, K\"{u}hn, Methuku and Osthus~\cite{KKKMO}.
In~\cite{MubayiV2}, we consider the analog of this for $k > 2$: let $\chi(e,k,q)$ denote the maximum possible chromatic number
of the hypergraph $H$ over all families $F$ of $e$ sets of size $q$ which pairwise intersect in at most one vertex. In~\cite{MubayiV2}, we determine the order of magnitude
of $\chi(e,k,q)$ up to a factor $q^{-2/(k - 1)}$ as $e \rightarrow \infty$:

\begin{thm} \label{coloring}
	Let $k,e \ge 3$ and $q \geq 1$. Then there exist absolute positive constants $c_1,c_2$  such that
 $$c_1\, \left(\frac{e^{1/2}}{\log e}\right)^{\frac{1}{k-1}} \, q^{1-\frac{2}{k-1}}  <
 \chi(e,k,q) < c_2\, \left(\frac{e^{1/2}}{\log e}\right)^{\frac{1}{k-1}}
 \,q$$
 where the lower bound holds for $e>q^{2k+3}$.
\end{thm}

It seems likely that $\chi(e,k,q) = \Theta(e^{1/2}/\log e) \cdot q$
as $e \rightarrow \infty$.

\end{document}